\documentclass[11pt]{amsart}
\usepackage{amsfonts}
\usepackage{setspace}
\usepackage{array}
\usepackage{booktabs}
\usepackage{diagbox} 
\usepackage{footnote}
\usepackage{url}
\usepackage{amssymb}
\usepackage{stmaryrd} 
\usepackage{mathrsfs}
\usepackage{latexsym,bm}
\usepackage{yfonts}
\usepackage{amsmath,amssymb}
\usepackage{geometry}
\usepackage{cite}
\usepackage{needspace}
\usepackage{hyperref}
\hypersetup{
	colorlinks=true,
	linkcolor=blue,
	citecolor=blue,
	urlcolor=blue,
	pdfauthor={},
	pdftitle={Schatten and weak Schatten classes of composition operator differences on Bergman spaces},
	pdfsubject={Composition operators on weighted Bergman spaces},
}
\newcommand{\declarationheading}[1]{%
  \par\addvspace{0.7\baselineskip}%
  \noindent\textbf{#1}\par\nobreak
  \vspace{0.15\baselineskip}\noindent\ignorespaces
}

\newtheorem{theorem}{Theorem}[section]
\newtheorem{lemma}[theorem]{Lemma}
\newtheorem{proposition}[theorem]{Proposition}

\theoremstyle{definition}

\newcommand{\D}{\mathbb D}
\newcommand{\dd}{\mathrm d}
\newcommand{\Id}{\operatorname{Id}}
\newcommand{\Toep}{\mathcal T}
\newcommand{\DP}{\mathcal D}
\DeclareMathOperator{\rank}{rank}
\DeclareMathOperator{\Tr}{Tr}
\DeclareMathOperator{\spanop}{span}
\catcode`@=11 \@addtoreset{equation}{section} \catcode`@=12
\title[Schatten and weak Schatten classes of composition operator differences]%
{Schatten and weak Schatten classes of composition operator differences on Bergman spaces}

\author{Qian Ai$^{1}$ }
\author{Rong Yang$^{*}$ }

\thanks{$^{1}$ {\small School of Mathematical Sciences, Soochow University, Suzhou 215006, P. R. China}}
\thanks{E-mail: 20234007001@stu.suda.edu.cn}

\thanks{$^{\rm *}$ {\small School of Mathematical Sciences, Dalian University of Technology, China}}
\thanks{E-mail: YRong@mail.dlut.edu.cn}

\numberwithin{equation}{section}

\begin{document}

\begin{abstract}
We investigate Schatten and weak Schatten class membership of differences of composition operators on the standard weighted Bergman spaces over the unit disk.
We prove that the known sufficient condition for Schatten $p$-class membership is also necessary for $0<p<2$, thereby completing the characterization for all $p>0$.
We also characterize when differences of composition operators belong to weak Schatten classes for every $p>0$.
Both results are expressed in terms of local-mass conditions involving the pseudohyperbolic distance between the inducing maps, and yield two-sided (quasi-)norm estimates.
\end{abstract}
	
	\maketitle

	
\section{Introduction and main results}
\label{sec:intro}

Let $\D$ be the unit disk in the complex plane $\mathbb C$, and let $H(\D)$ denote the space of holomorphic functions on $\D$.
A holomorphic self-map $\varphi$ of $\D$ defines the composition operator $C_\varphi$ on $H(\D)$ by
\[
C_\varphi f=f\circ\varphi,\qquad f\in H(\D).
\]
Much of the theory concerns how the geometry and boundary behavior of $\varphi$ determine the boundedness, compactness, and spectral properties of $C_\varphi$ on spaces of holomorphic functions.
The monographs of Cowen and MacCluer \cite{CowenMacCluer1995} and Shapiro \cite{Shapiro1993} provide systematic accounts of this connection between operator theory and classical function theory.

One topic of particular interest in the theory of composition operators is the study of differences of two composition operators and, more generally, linear combinations of composition operators.
This line of research grew out of investigations of the operator norm topology on the space of composition operators, initiated in the Hardy space setting by Berkson \cite{Berkson} and developed further by Shapiro and Sundberg \cite{Shapiro}.
Two distinct noncompact composition operators may have a compact difference, so the interaction between their symbols is essential to the compactness problem.
On standard weighted Bergman spaces, compact differences have been characterized by boundary conditions involving the pseudohyperbolic distance between the symbols \cite{Moorhouse2005}.
Carleson measure methods have also yielded boundedness and compactness criteria for differences between standard weighted Bergman spaces and on weighted Bergman spaces over the unit ball, as well as a compactness characterization on the Hardy space; see \cite{Saukko2011,KooWang2014,ChoeChoiKooYang2020}.

Schatten and weak Schatten classes play an important role in operator theory on Hilbert spaces, providing quantitative information about singular value decay, whereas compactness alone gives no rate of decay.
Gohberg and Krein \cite{GohbergKrein1969} and Simon \cite{Simon2005} developed the theory of operator ideals on Hilbert spaces, while Pietsch \cite{Pietsch1987} studied eigenvalues and $s$-numbers and Birman and Solomyak \cite{BirmanSolomyak1977} established strong and weak Schatten estimates for integral operators.
For differences of composition operators on standard weighted Bergman spaces, Choe, Hosokawa and Koo \cite{ChoeHosokawaKoo2011} characterized the Hilbert--Schmidt property.
Using local averages of a joint pullback measure, Choe, Choi, Koo and Park \cite{ChoeChoiKooPark2026} studied Schatten $p$-class differences for $p\geq2$, and proved that the same condition is sufficient for $0<p<2$.
They left open whether this condition is also necessary for $0<p<2$ \cite{ChoeChoiKooPark2026}.
Motivated by this question, we prove the missing necessity and give both strong and weak Schatten membership for every $p>0$ in terms of the same local averages.

For individual composition operators, Luecking and Zhu \cite{LueckingZhu1992} obtained Schatten class criteria on Hardy and weighted Bergman spaces.
Bayart, Gilmore and \c{S}ahin \cite{BayartGilmoreSahin2026} established upper and lower estimates for approximation numbers of differences of composition operators on the Hardy space.
For differences of weighted composition operators, Du, Tong and Yang \cite[Theorem~1.3]{DuTongYang2026} gave related local-average conditions, with sufficiency for all $p>0$ and necessity for $p\geq2$ under their stated parameter assumptions.
Related Schatten class results on weighted Bergman spaces over the unit ball were proved by Hu and Yang \cite{HuYang2026}.

To begin with, we introduce the function spaces and operator ideals involved.
For $\alpha>-1$, put
\[
\dd A_\alpha(z)=(\alpha+1)(1-|z|^2)^\alpha\dd A(z),
\qquad \dd A(z)=\pi^{-1}\dd x\,\dd y.
\]
The standard weighted Bergman space $A_\alpha^2$ consists of all $f\in H(\D)$ for which
\[
\|f\|_{A_\alpha^2}^2
=\int_\D |f(z)|^2\,\dd A_\alpha(z)<\infty.
\]
With this norm, $A_\alpha^2$ is a Hilbert space, and Littlewood's subordination principle ensures that every $C_\varphi$ is bounded on $A_\alpha^2$; see \cite{CowenMacCluer1995,DurenSchuster2004}.
For a compact operator $T$ between Hilbert spaces, let $s_n(T)$ denote its singular values in nonincreasing order.
For $p>0$, set
\[
\|T\|_{S_p}^p=\sum_{n\geq1}s_n(T)^p,
\qquad
\|T\|_{S_{p,\infty}}=\sup_{n\geq1}n^{1/p}s_n(T).
\]
The Schatten class $S_p$ and weak Schatten class $S_{p,\infty}$ consist of the compact operators for which the corresponding quantities are finite, respectively.
In particular, $S_2$ is the Hilbert--Schmidt class.

To state our results, we introduce some notation.
Let $\varphi,\psi$ be holomorphic self-maps of $\D$, and write $D=C_\varphi-C_\psi$.
For $a,b\in\D$ and $z\in\D$, put
\[
\rho(a,b)=\left|\frac{a-b}{1-\overline a b}\right|,
\qquad \rho_z=\rho(\varphi(z),\psi(z)).
\]
Define the joint pullback measure by
\begin{equation}\label{eq:joint-measure}
\mu(E)=\int_\D\rho_z^2
\bigl(\mathbf1_E(\varphi(z))+\mathbf1_E(\psi(z))\bigr)\,\dd A_\alpha(z),
\end{equation}
where $E\subset\D$ is Borel and $\mathbf1_E$ denotes its characteristic function.
Then $\mu$ is a positive finite Borel measure with $\mu(\D)\leq2$.
Let $d(a,b)=\operatorname{arctanh}\rho(a,b)$ and $B(a,r)=\{b\in\D:d(a,b)<r\}$ for $r>0$.
Fix $r>0$ and let $\{a_j\}$ be a lattice of mesh $r$, that is, a maximal $r$-separated set in the metric $d$.
Define the local densities by
\begin{equation}\label{eq:lattice-density}
m_j=\frac{\mu(B(a_j,r))}{(1-|a_j|^2)^{\alpha+2}}.
\end{equation}
Since the weighted area of $B(a,r)$ is comparable to $(1-|a|^2)^{\alpha+2}$ for fixed $r$, these densities are comparable to the local averages used in \cite{ChoeChoiKooPark2026}.
We use the usual sequence (quasi-)norms $\|\cdot\|_{\ell^p}$ and $\|\cdot\|_{\ell^{p,\infty}}$, with
\[
\|x\|_{\ell^{p,\infty}}
=\sup_{t>0}t\,\#\{j:|x_j|>t\}^{1/p}.
\]

Throughout the paper, $C$ denotes a positive constant that may change from one occurrence to another.
For nonnegative quantities $X$ and $Y$, we write $X\lesssim Y$ or $Y\gtrsim X$ if $X\leq CY$ for some constant $C>0$, and $X\asymp Y$ if both $X\lesssim Y$ and $Y\lesssim X$.

The following are our main results. Our first theorem establishes the missing necessity for $0<p<2$.

\begin{theorem}\label{thm:strong}
Let $\alpha>-1$, $0<p<2$, and $\varphi,\psi$ be holomorphic self-maps of $\D$.
For $r>0$, let $\{a_j\}$ be a $r$-lattice.
Then $C_\varphi-C_\psi\in S_p(A_\alpha^2)$ if and only if $(m_j^{1/2})\in \ell^p$.
In this case,
\begin{equation}\label{eq:main-strong-norm}
\|C_\varphi-C_\psi\|_{S_p}
\asymp\|(m_j^{1/2})\|_{\ell^{p}},
\end{equation}
where the comparison constants depend only on $\alpha$, $p$, and $r$.
\end{theorem}

Together with \cite[Theorem~1.1]{ChoeChoiKooPark2026}, Theorem~\ref{thm:strong} completes the strong Schatten class characterization for every $p>0$.
Combining localization, derivative interpolation for nearby image points, sparse-matrix estimates, and uniform reconstruction from local subspaces spanned by reproducing kernels, we also establish weak Schatten membership of differences of composition operators for every $p>0$.

\begin{theorem}\label{thm:weak}
Let $\alpha>-1$, $p>0$, and $\varphi,\psi$ be holomorphic self-maps of $\D$. For $r>0$, let $\{a_j\}$ be a $r$-lattice. Then $C_\varphi-C_\psi\in S_{p,\infty}(A_\alpha^2)$ if and only if $(m_j^{1/2})\in\ell^{p,\infty}$.
In this case,
\begin{equation}\label{eq:main-weak-norm}
\|C_\varphi-C_\psi\|_{S_{p,\infty}}
\asymp\|(m_j^{1/2})\|_{\ell^{p,\infty}},
\end{equation}
where the comparison constants depend only on $\alpha$, $p$, and $r$.
\end{theorem}

The paper is organized as follows. Section~\ref{sec:local-tests} develops the kernel and lattice estimates and derives necessary bounds by local testing.
Section~\ref{sec:tools} establishes the sparse matrix estimates and the uniform reconstruction of local kernel subspaces.
Section~\ref{sec:necessity} proves the remaining necessary estimates by treating close and separated image points.
Section~\ref{sec:toeplitz} gives both sufficient directions and completes the proofs of the main theorems.

\section{Kernels and local tests}
\label{sec:local-tests}

In this section, we establish estimates for reproducing kernels and compare the local density conditions associated with different lattices.
These estimates will be used to obtain necessary bounds for Schatten and weak Schatten membership of differences of composition operators.
For standard estimates on weighted Bergman spaces, see \cite[Chapter~4]{Zhu2007}.

\subsection{Kernel estimates and lattices}

Let $H=A_\alpha^2$.
For convenience, we use the following notation throughout the rest of the paper:
\[
\beta=\alpha+2>1,\qquad \Delta(a)=1-|a|^2.
\]
For $a\in\D$, the reproducing kernel $K_a$, its norm $s_a$, and the normalized kernel $k_a$ are given by
\[
K_a(z)=(1-\overline a z)^{-\beta},\qquad
s_a=\|K_a\|=\Delta(a)^{-\beta/2},\qquad k_a=s_a^{-1}K_a.
\]
The reproducing identity is $f(a)=\langle f,K_a\rangle$ for $f\in H$.

For local testing, we use the operator $\DP=JD:H\to L^2(\dd A_\alpha)$, where $J$ is the natural isometric inclusion.
Thus $\DP f=f\circ\varphi-f\circ\psi$ for $f\in H$.
For $u,v\in H$, let $u\otimes v$ denote the operator on $H$ defined by $(u\otimes v)f=\langle f,v\rangle u$.

\begin{lemma}
\label{lem:analysis-operator}
One has $\DP^*\DP=D^*D$, so $\DP$ and $D$ have the same singular values whenever either is compact.
Furthermore,
\begin{equation}
\label{eq:pair-energy}
D^*D=\int_\D (K_{\varphi(z)}-K_{\psi(z)})
\otimes(K_{\varphi(z)}-K_{\psi(z)})\,\dd A_\alpha(z)
\end{equation}
in the weak operator sense.
\end{lemma}

\begin{proof}
The first identity is $D^*J^*JD=D^*D$.
For $f,g\in H$, the form on the right of \eqref{eq:pair-energy} is $\int_\D(\DP f)\overline{\DP g}\,\dd A_\alpha$.
It is absolutely convergent by Cauchy--Schwarz.
This proves the second identity without requiring norm convergence of the operator-valued integral.
\end{proof}

We shall use the comparisons
\begin{equation}
\label{eq:local-comparability}
\Delta(a)\asymp_r\Delta(c),\qquad
|1-\overline a b|\asymp_r|1-\overline c b|\quad(b\in\D),
\end{equation}
valid whenever $d(a,c)\leq r$.
Both follow from the disk automorphism $\sigma_c(z)=(c-z)/(1-\overline c z)$.

We recall the following kernel estimate from \cite[Proposition~3.5]{ChoeHosokawaKoo2011}.
We include a proof for completeness.

\begin{lemma}
\label{lem:kernel-difference}
For all $a,b\in\D$,
\[
\|K_a-K_b\|^2\asymp_\alpha
\rho(a,b)^2(s_a^2+s_b^2).
\]
\end{lemma}

\begin{proof}
The identities
\[
\|K_a-K_b\|^2=s_a^2+s_b^2-2\operatorname{Re}(1-\overline b a)^{-\beta},
\qquad
|\langle k_a,k_b\rangle|=(1-\rho(a,b)^2)^{\beta/2}
\]
give the lower bound when $\rho(a,b)\geq1/2$, since $|\langle k_a,k_b\rangle|\leq(3/4)^{\beta/2}<1$ and $2s_as_b\leq s_a^2+s_b^2$.
The upper bound follows from the triangle inequality and $\rho(a,b)^2\geq1/4$.

Suppose now that $\rho(a,b)<1/2$.
By \eqref{eq:local-comparability}, $\Delta(a)\asymp\Delta(b)$, hence $s_a\asymp s_b$.
The lower bound follows from
\[
\begin{aligned}
\|K_a-K_b\|^2
&\geq s_a^2+s_b^2-2s_as_b(1-\rho(a,b)^2)^{\beta/2}\\
&\geq2s_as_b\bigl[1-(1-\rho(a,b)^2)^{\beta/2}\bigr]
\gtrsim \rho(a,b)^2s_a^2.
\end{aligned}
\]
For the upper bound, write $K_a^{[1]}=\partial_{\overline a}K_a$.
Direct differentiation gives
\begin{equation}
\label{eq:derivative-kernel-norm}
K_a^{[1]}(z)=\frac{\beta z}{(1-\overline a z)^{\beta+1}},\qquad
\|K_a^{[1]}\|^2=
\beta(1+\beta|a|^2)\Delta(a)^{-\beta-2}.
\end{equation}
The pseudohyperbolic disk $\{u:\rho(a,u)<1/2\}$ is a Euclidean disk containing the line segment from $a$ to $b$, and $\Delta(u)\asymp\Delta(a)$ throughout this disk.
Moreover, $|a-b|\lesssim\rho(a,b)\Delta(a)$.
Integrating the derivative of the kernel map along this segment yields $\|K_a-K_b\|\lesssim\rho(a,b)s_a$, and the upper bound follows.
\end{proof}

Lemma~\ref{lem:kernel-difference}, together with Tonelli's theorem and Parseval's identity, yields the Hilbert--Schmidt estimate of \cite[Theorem~3.6]{ChoeHosokawaKoo2011}:
\[
\|D\|_{S_2}^2
=\int_\D\|K_{\varphi(z)}-K_{\psi(z)}\|^2\,\dd A_\alpha(z)
\asymp\int_\D s_a^2\,\dd\mu(a),
\]
with both sides allowed to be infinite.

For a lattice of mesh $\eta>0$, choose nearest-center cells $Q_j$, assigning every tie to the least index.
Then the cells are Borel, form an exact disjoint partition, and
\[
B(a_j,\eta/2)\subset Q_j\subset B(a_j,\eta),\qquad
\int_{Q_j}\dd A_\alpha\asymp_{\alpha,\eta}\Delta(a_j)^\beta.
\]
The tie convention is fixed on all points, not merely almost everywhere, since pushforward measures may have atoms.
Put
\[
s_j=\Delta(a_j)^{-\beta/2},\qquad
q_j=s_j^2\mu(Q_j)=q_{j,\varphi}+q_{j,\psi},
\]
where
\[
q_{j,\varphi}=s_j^2\int_{\varphi^{-1}(Q_j)}\rho_z^2\,\dd A_\alpha(z),
\qquad
q_{j,\psi}=s_j^2\int_{\psi^{-1}(Q_j)}\rho_z^2\,\dd A_\alpha(z).
\]

\begin{lemma}
\label{lem:discretization}
Fix two lattices with arbitrary fixed positive meshes.
Let $(q_j)$ denote the cell densities for the first lattice and $(m_j)$ the disk densities for the second; the index on each side refers to its own lattice.
Then, for every $p>0$,
\begin{align}
\|(q_j^{1/2})\|_{\ell^p}&\asymp\|(m_j^{1/2})\|_{\ell^p},
\label{eq:strong-discretization}\\
\|(q_j^{1/2})\|_{\ell^{p,\infty}}&\asymp
\|(m_j^{1/2})\|_{\ell^{p,\infty}}.
\label{eq:weak-discretization}
\end{align}
The constants depend only on $\alpha$, the two fixed meshes, and $p$.
Moreover, for every fixed $R>0$, a lattice of mesh $\eta$ is a union of finitely many $R$-separated sequences, with the number of colors bounded in terms of $R$ and $\eta$ alone.
\end{lemma}

\begin{proof}
A disk of one fixed radius meets only boundedly many cells of a fixed mesh, their centers have comparable depths, and the incidence multiplicity is bounded.
Thus each density is bounded by a constant times a sum of boundedly many densities of the other kind, with bounded reverse multiplicity.
The same statement for two different cell partitions follows by intersecting the cells.
Raising such finite sums to the power $p/2$ and summing proves \eqref{eq:strong-discretization}.
If a finite sum of nonnegative terms exceeds $t^2$, at least one term exceeds a fixed multiple of $t^2$.
Counting the corresponding indices proves \eqref{eq:weak-discretization}.

For the last assertion, join two centers whose distance is less than $R$.
Disjoint $\eta/2$-balls and invariant area give a uniform bound for the degree of this graph.
Indeed, with $\dd\lambda(z)=\Delta(z)^{-2}\dd A(z)$, $\lambda(B(a,t))=\sinh^2t$.
Greedy coloring of the countable graph therefore uses only finitely many colors.
Centers of the same color are $R$-separated.
\end{proof}

\subsection{Local testing}

We use the normalized kernels and first-derivative kernels to estimate the cell densities $q_{j,\varphi}$ and $q_{j,\psi}$.
Restricting their images under $\DP$ to disjoint sets gives the strong necessary estimate for $p\geq1$ and the weak necessary estimate for $p>1$.

For fixed $r>0$, $a\in\D$, and $f\in H$, local submean and Cauchy estimates give
\begin{equation}
\label{eq:local-derivatives}
|f^{(m)}(a)|^2\leq C_{\alpha,m,r}\Delta(a)^{-\beta-2m}
\int_{B(a,r)}|f|^2\,\dd A_\alpha,\qquad m\geq0.
\end{equation}
To obtain this estimate, apply the submean and Cauchy estimates on two concentric Euclidean disks centered at $a$, with radii comparable to $\Delta(a)$ and contained in $B(a,r)$.
The weight $(1-|z|^2)^\alpha$ is comparable to $\Delta(a)^\alpha$ on these disks.

Let $k_c^{[0]}=k_c$ and $k_c^{[1]}=K_c^{[1]}/\|K_c^{[1]}\|$.
The next lemma gives the pointwise estimate and Bessel bounds needed for these tests.

\begin{lemma}
\label{lem:two-kernel-test}
For fixed $r>0$, if $a\in B(c,r)$ and $b\in\D$, then
\begin{equation}
\label{eq:two-kernel-test}
\sum_{v=0}^1|k_c^{[v]}(a)-k_c^{[v]}(b)|^2
\geq c_{\alpha,r}\rho(a,b)^2s_c^2.
\end{equation}
Moreover, for every $\sigma>0$, each of the two families indexed by a $\sigma$-separated sequence is Bessel, with bound depending only on $\alpha$ and $\sigma$.
The same bound, up to an absolute factor, holds for arbitrary choices between the two families with at most two vectors associated with each center.
\end{lemma}

\begin{proof}
Set
\begin{align*}
A&=1-\overline c a, & B&=1-\overline c b,
& \tau&=\frac AB, & \chi&=A^\beta B^{-\beta},\\
x&=1-\chi, & y&=a-b\tau\chi.
\end{align*}
Since $\operatorname{Re}A,\operatorname{Re}B>0$, the powers are taken with the holomorphic branch of the logarithm on the right half-plane.
A direct calculation gives
\[
y-b\tau x=\frac{a-b}{1-\overline c b}.
\]
By \eqref{eq:local-comparability},
\[
|A|\asymp_r\Delta(c),\qquad
|B|\asymp_r|1-\overline a b|,\qquad |\tau|\leq C_r.
\]
Hence
\[
\rho(a,b)
\lesssim_r |y|+|x|,
\]
and therefore $|x|^2+|y|^2\gtrsim_r\rho(a,b)^2$.
Moreover,
\begin{align*}
k_c(a)-k_c(b)
&=\Delta(c)^{\beta/2}A^{-\beta}x,\\
k_c^{[1]}(a)-k_c^{[1]}(b)
&=\frac{\sqrt\beta\,\Delta(c)^{(\beta+2)/2}}
{\sqrt{1+\beta|c|^2}}A^{-\beta-1}y.
\end{align*}
Thus the two differences are comparable in modulus to $s_c|x|$ and $s_c|y|$, respectively, which proves \eqref{eq:two-kernel-test}.

For the Bessel assertion, let $\{a_j\}$ be $\sigma$-separated.
For $v=0,1$, the reproducing property and \eqref{eq:local-derivatives}, on a fixed ball radius depending only on $\sigma$, give
\[
|\langle f,k_{a_j}^{[v]}\rangle|^2
\leq C_{\alpha,\sigma}
\int_{B(a_j,r_1)}|f|^2\,\dd A_\alpha.
\]
The balls have uniformly bounded overlap, so summation over $j$ gives a Bessel bound.
Taking a subfamily and allowing at most two repetitions at each center changes the bound only by a fixed factor.
\end{proof}

The local tests yield the following necessary estimates.

\begin{proposition}
\label{prop:banach-necessity}
For every lattice of mesh $\eta>0$, if $D\in S_p$, then
\[
\sum_jq_j^{p/2}\leq C_{\alpha,p,\eta}\|D\|_{S_p}^p,
\qquad 1\leq p<\infty.
\]
Moreover, if $D\in S_{p,\infty}$, then
\begin{equation}
\label{eq:weak-necessity-large-p}
\|(q_j^{1/2})\|_{\ell^{p,\infty}}
\leq C_{\alpha,p,\eta}\|D\|_{S_{p,\infty}},
\qquad p>1.
\end{equation}
\end{proposition}

\begin{proof}
For $i=(j,\xi)$, where $\xi\in\{\varphi,\psi\}$, put $\Omega_i=\xi^{-1}(Q_j)$ and $q_i=q_{j,\xi}$.
Integrating \eqref{eq:two-kernel-test} on $\Omega_i$ and choosing the larger of its two terms produces a unit vector $F_i\in\{k_{a_j}^{[0]},k_{a_j}^{[1]}\}$ such that
\begin{equation}
\label{eq:selected-tests}
d_i:=\|\mathbf1_{\Omega_i}\DP F_i\|_2\geq c q_i^{1/2},
\qquad
\left\|\sum_i z_iF_i\right\|\leq B\|(z_i)\|_2.
\end{equation}
The constants depend only on $\alpha$ and $\eta$.
Indices with $d_i=0$ also have $q_i=0$ and are omitted.

Use the Hilbert space $\mathcal K=L^2(\dd A_\alpha)\oplus L^2(\dd A_\alpha)$ and set $\widetilde\DP f=2^{-1/2}(\DP f,\DP f)$.
Its singular values equal those of $D$ whenever $D$ is compact.
Put $u_i=\mathbf1_{\Omega_i}\DP F_i/d_i$ in the first output copy for $\xi=\varphi$ and the second for $\xi=\psi$.
These vectors are orthonormal, since the sets $\Omega_i$ in each copy are disjoint, and
\[
\langle\widetilde\DP F_i,u_i\rangle=2^{-1/2}d_i.
\]
For a finite index set $E$, let $A_Ee_i=F_i$ and $U_Ee_i=u_i$.
Then $\|A_E\|\leq B$ and $U_E$ is an isometry.
This construction requires only boundedness of $D$ and will be used for both estimates.

First suppose that $1\leq p<\infty$ and $D\in S_p$.
Set $V_{E,w}=A_E\operatorname{diag}(w_i)U_E^*$.
The ideal property gives
\[
\|V_{E,w}\|_{S_{p'}}\leq B\|w\|_{\ell^{p'}},
\]
where $p'$ is the H\"older conjugate of $p$, with $p'=\infty$ when $p=1$; $S_\infty$ denotes the operator norm.
As $V_{E,w}$ has finite rank, cyclicity of the finite-rank trace yields
\[
\Tr(V_{E,w}\widetilde\DP)
=\Tr\!\left(\operatorname{diag}(w_i)U_E^*\widetilde\DP A_E\right)
=2^{-1/2}\sum_{i\in E}w_id_i.
\]
Schatten H\"older and finite dimensional sequence duality imply
\[
\left(\sum_{i\in E}d_i^p\right)^{1/p}\leq\sqrt2 B\|D\|_{S_p}.
\]
This includes $p=1$ by using the operator norm for $V_{E,w}$.
Let $E$ increase and use \eqref{eq:selected-tests}.
Finally,
\[
(q_{j,\varphi}+q_{j,\psi})^{p/2}
\leq 2^{\max\{p/2-1,0\}}
\bigl(q_{j,\varphi}^{p/2}+q_{j,\psi}^{p/2}\bigr),
\]
which proves the strong estimate.

Now suppose that $p>1$ and $D\in S_{p,\infty}$.
For a finite index set $E$ of cardinality $N$, put $V_E=A_EU_E^*$.
Then $\rank(V_E)\leq N$ and $\|V_E\|\leq B$.
The finite-rank trace inequality gives
\[
\begin{aligned}
2^{-1/2}\sum_{i\in E}d_i
&=\Tr(V_E\widetilde\DP)
\leq\sum_{n=1}^N s_n(V_E)s_n(D)\\
&\leq B\|D\|_{S_{p,\infty}}\sum_{n=1}^N n^{-1/p}
\leq C_p B N^{1-1/p}\|D\|_{S_{p,\infty}}.
\end{aligned}
\]
Taking $E$ among the indices for which $d_i>t$, then taking the supremum over finite $E$ and $t>0$, yields $\|(d_i)\|_{\ell^{p,\infty}}\leq C_pB\|D\|_{S_{p,\infty}}$.
By \eqref{eq:selected-tests}, $q_i^{1/2}\lesssim d_i$.
Combining the two sequences $(q_{j,\varphi}^{1/2})$ and $(q_{j,\psi}^{1/2})$ by a level-set union bound proves \eqref{eq:weak-necessity-large-p}.
\end{proof}

\Needspace{9\baselineskip}
\section{Sparse matrices and reconstruction}
\label{sec:tools}

We establish Schatten estimates for sparse matrices and construct bounded reconstruction operators for subspaces spanned by reproducing kernels.

For a sequence $x=(x_j)$, let $x^*=(x_n^*)$ denote the nonincreasing rearrangement of $(|x_j|)$.

\subsection{Sparse matrix estimates}

The localization argument below produces finite matrices with a uniformly bounded number of nonzero entries in each row.
We record strong and weak estimates whose constants are independent of the matrix dimensions.

\begin{lemma}
\label{lem:sparse-strong}
Let $m\geq1$ be an integer, and let $B:\mathbb C^N\to\mathbb C^M$ have at most $m$ nonzero entries in each row.
Put $h_j=\|Be_j\|^2$.
If $0<v<2$, then
\begin{equation}
\label{eq:sparse-strong}
m^{v/2-1}\sum_jh_j^{v/2}
\leq\|B\|_{S_v}^v\leq\sum_jh_j^{v/2}.
\end{equation}
Moreover,
\begin{equation}
\label{eq:sparse-pointwise-upper}
s_n(B)\leq\sqrt m\,(h^{1/2})_n^*,\qquad n\geq1.
\end{equation}
Both constants in \eqref{eq:sparse-strong} are optimal uniformly over this matrix class.
\end{lemma}

\begin{proof}
Delete the zero columns and let $H_0=\operatorname{diag}(h_j)$.
Write
\[
B^*B=H_0^{1/2}QH_0^{1/2},
\qquad
Q=H_0^{-1/2}B^*BH_0^{-1/2}.
\]
Then $Q\geq0$ and $Q_{jj}=1$.
If $B=(b_{ij})$, then for $x\in\mathbb C^N$,
\begin{align*}
\langle Qx,x\rangle
&=\sum_i\left|\sum_{j:b_{ij}\neq0}
\frac{b_{ij}}{\sqrt{h_j}}x_j\right|^2\\
&\leq m\sum_i\sum_j\frac{|b_{ij}|^2}{h_j}|x_j|^2
=m\|x\|_2^2.
\end{align*}
Thus $0\leq Q\leq mI$.

Put $\theta=v/2\in(0,1)$.
In the notation of \cite[Theorem~1]{Audenaert2008}, apply the Araki--Lieb--Thirring inequality with $A=H_0^{1/2}$, $B=Q$, $r=\theta$, and $q=1$.
This gives
\[
\Tr\!\left[(H_0^{1/2}QH_0^{1/2})^\theta\right]
\geq \Tr(H_0^\theta Q^\theta).
\]
Since $0\leq Q\leq mI$ and $t^\theta\geq m^{\theta-1}t$ on $[0,m]$,
\[
\Tr\!\left[(B^*B)^\theta\right]
\geq m^{\theta-1}\Tr(H_0^\theta Q)
=m^{\theta-1}\sum_jh_j^\theta.
\]
For the reverse inequality, let $B^*B=\sum_k\lambda_k u_k\otimes u_k$.
Scalar concavity of $t^\theta$ gives
\[
\langle (B^*B)^\theta e_j,e_j\rangle
\leq \langle B^*Be_j,e_j\rangle^\theta=h_j^\theta.
\]
Summation over $j$ yields the upper bound in \eqref{eq:sparse-strong}.
Finally, $Q\leq mI$ implies $B^*B\leq mH_0$, and the min--max principle gives \eqref{eq:sparse-pointwise-upper}.
A single row containing $m$ ones attains the lower constant, while a diagonal matrix attains the upper one.
\end{proof}

The corresponding weak estimate will be applied only to finite matrices, so it requires no a priori summability.

\begin{lemma}
\label{lem:sparse-weak}
For every $p>0$, matrices as in Lemma~\ref{lem:sparse-strong} satisfy
\[
C_{p,m}^{-1}\|(h_j^{1/2})\|_{\ell^{p,\infty}}
\leq \|B\|_{S_{p,\infty}}
\leq \sqrt m\,\|(h_j^{1/2})\|_{\ell^{p,\infty}}.
\]
\end{lemma}

\begin{proof}
The upper estimate follows immediately from \eqref{eq:sparse-pointwise-upper}.
For the reverse estimate, fix $0<v<\min\{p,2\}$ and let $J\subset\{j:h_j^{1/2}>t\}$ be finite, with $|J|=N_1$.
The column restriction $B_J$ still has row support at most $m$, hence Lemma~\ref{lem:sparse-strong} gives
\[
m^{v/2-1}N_1t^v
\leq\|B_J\|_{S_v}^v
\leq\sum_{n=1}^{N_1}s_n(B)^v
\leq C_{p,v}\|B\|_{S_{p,\infty}}^vN_1^{1-v/p}.
\]
Therefore
\[
tN_1^{1/p}\leq C_{p,m}\|B\|_{S_{p,\infty}}.
\]
Taking the supremum over finite subsets of the level set, and then over $t>0$, proves the lower estimate.
The restriction here is purely a finite matrix column restriction and has no analytic interpolation content.
\end{proof}

\subsection{Uniform reconstruction}

We next construct a uniformly bounded reconstruction operator for subspaces spanned by one or two reproducing kernels.
The geometric situation is related to generalized interpolation without separation \cite{Luecking2014}, but the estimate needed here is proved directly from kernel identities and a weighted Schur test.
In particular, its constant is independent of the distance between two points in the same local cluster.

We begin with a weighted estimate for the off-diagonal kernel terms.
Put
\[
\kappa_\beta(a,b)=|\langle k_a,k_b\rangle|
=\frac{\Delta(a)^{\beta/2}\Delta(b)^{\beta/2}}
{|1-\overline a b|^\beta},\qquad w(a)=\Delta(a)^{1/2}.
\]

\begin{lemma}
\label{lem:weighted-tail}
Fix $\eta_0>0$.
Let $(b_j)$ be a sequence in which each point occurs at most twice and any two distinct points are separated by at least $\eta_0$.
There is a constant $C_{\beta,\eta_0}$ such that, for $L\geq1$,
\[
\sup_{a\in\D}\sum_{d(a,b_j)\geq L}
\kappa_\beta(a,b_j)\frac{w(b_j)}{w(a)}
\leq C_{\beta,\eta_0}(1+L)e^{-(\beta-1)L}.
\]
\end{lemma}

\begin{proof}
Small fixed hyperbolic disks about the distinct $b_j$ are disjoint.
On each such disk, the summand, with $b_j$ replaced by the integration variable, is comparable to its value at $b_j$, by \eqref{eq:local-comparability}.
Invariant area of each disk is a fixed positive constant.
It is consequently enough, with a fixed reduction of $L$, to estimate
\[
\int_{d(a,z)\geq L}
\kappa_\beta(a,z)\frac{w(z)}{w(a)}\,\dd\lambda(z).
\]
Set $z=\sigma_a(u)$.
Invariance of $\dd\lambda$ and $\Delta(\sigma_a(u))=\Delta(a)\Delta(u)/|1-\overline a u|^2$ turn this integral into
\begin{equation}
\label{eq:tail-integral}
\int_{|u|\geq\tanh L}
\frac{\Delta(u)^{(\beta-3)/2}}{|1-\overline a u|}\,\dd A(u).
\end{equation}
Uniformly in $|a|\leq1$, angular integration at radius $t$ is at most $C\log(e/(1-t))$.
To see this, split the angle at $1-|a|t$ and use $|1-|a|te^{i\theta}|\gtrsim (1-|a|t)+|\theta|$ for $|\theta|\leq\pi$, treating $|a|t\leq1/2$ separately.
With $v=1-t^2$, \eqref{eq:tail-integral} is therefore bounded by
\[
C\int_0^{\operatorname{sech}^2L}
v^{(\beta-3)/2}\log(e/v)\,\dd v
\leq C_\beta(1+L)e^{-(\beta-1)L}.
\]
The exponent is integrable precisely because $\beta>1$.
The fixed reduction of $L$ and the possible multiplicity two change only the constant.
\end{proof}

For nearby kernel points, we use orthonormal coordinates in the span of their reproducing kernels.
For $a\in\D$, define
\[
(U_af)(z)=\frac{\Delta(a)^{\beta/2}}{(1-\overline a z)^\beta}
f(\sigma_a(z)).
\]
Fix $\tau_0=1/8$ and, for $0<|u|\leq\tau_0$, put
\[
v_{a,u}=U_a\frac{K_u-1}{\|K_u-1\|}.
\]

\begin{lemma}
\label{lem:close-atom-envelope}
For $a\in\D$ and $0<|u|\leq\tau_0$, the vectors $k_a,v_{a,u}$ form an orthonormal basis of $\spanop\{K_a,K_{\sigma_a(u)}\}$.
Associate the center $a$ to either of these vectors.
If $e$ and $f$ are such vectors with centers $a$ and $b$, respectively, then
\begin{equation}
\label{eq:close-envelope}
|\langle e,f\rangle|\leq C_\beta\kappa_\beta(a,b).
\end{equation}
The constant is independent of the distances between the two kernel points in either pair.
\end{lemma}

\begin{proof}
The automorphism change of variables in $\dd A_\alpha$ shows that $U_a$ is a unitary involution, with the analytic branch fixed at zero.
In particular, $U_a1=k_a$.
Since $1$ is orthogonal to $K_u-1$ and
\[
\|K_u-1\|^2=(1-|u|^2)^{-\beta}-1\asymp_\beta|u|^2,
\]
the vectors $k_a,v_{a,u}$ are orthonormal.
The identity
\[
U_aK_u
=\frac{\Delta(a)^{\beta/2}}{(1-a\overline u)^\beta}
K_{\sigma_a(u)}
\]
shows that they span the stated subspace.

Set
\[
F_{a,b}(u,v)=\langle U_aK_u,U_bK_v\rangle.
\]
For $|u|,|v|\leq2\tau_0$, the two kernel points remain in fixed hyperbolic neighborhoods of $a$ and $b$.
The preceding identity and \eqref{eq:local-comparability} give
\begin{equation}
\label{eq:translated-kernel-envelope}
|F_{a,b}(u,v)|\leq C_\beta\kappa_\beta(a,b).
\end{equation}
As a function of $\overline u$ and $v$, $F_{a,b}(u,v)$ is analytic on the bidisk $|u|,|v|<1$.
Applying Cauchy's estimate on the larger bidisk $|u|,|v|\leq2\tau_0$ to points with $|u|,|v|\leq\tau_0$ gives
\begin{align*}
|F_{a,b}(u,0)-F_{a,b}(0,0)|
&\leq C_\beta|u|\kappa_\beta(a,b),\\
|F_{a,b}(0,v)-F_{a,b}(0,0)|
&\leq C_\beta|v|\kappa_\beta(a,b),\\
|F_{a,b}(u,v)-F_{a,b}(u,0)-F_{a,b}(0,v)+F_{a,b}(0,0)|
&\leq C_\beta|u||v|\kappa_\beta(a,b).
\end{align*}
Now $k_a=U_a1$ and $v_{a,u}=U_a(K_u-1)/\|K_u-1\|$, with $\|K_u-1\|\asymp_\beta|u|$.
Estimate \eqref{eq:translated-kernel-envelope} handles $(k_a,k_b)$, the first two difference estimates handle $(v_{a,u},k_b)$ and $(k_a,v_{b,v})$, and the mixed difference handles $(v_{a,u},v_{b,v})$.
Division by the corresponding normalizing factors proves \eqref{eq:close-envelope}.
\end{proof}

The preceding estimates give the following reconstruction result.

\begin{proposition}
\label{prop:cluster-interpolation}
There are $R_0=R_0(\alpha)$ and $A_0=A_0(\alpha)$ with the following property.
Let $R\geq R_0$, and let $\Gamma_0,\Gamma_1$ be disjoint $R$-separated sequences.
Join two points in different sequences when their hyperbolic distance is less than $R/4$.
The components $C$ have at most two points.
Put
\[
H_C=\spanop\{K_a:a\in C\},\qquad
\mathcal X=\bigoplus_C H_C,
\]
where the direct sum has the original Hilbert norm on each $H_C$.
If $P_C:H\to H_C$ is orthogonal projection, the operator
\begin{equation}
\label{eq:cluster-restriction}
\mathcal R f=(P_Cf)_C:H\longrightarrow\mathcal X
\end{equation}
is bounded and onto.
It has a linear right inverse $\mathcal A$ such that
\[
\mathcal R\mathcal A=\Id,\qquad \|\mathcal A\|\leq A_0.
\]
The constants are independent of the sequences and of all distances between the two points of a component.
The same assertion holds for one $R$-separated sequence alone.
\end{proposition}

\begin{proof}
A point cannot have two neighbors: those neighbors would have the same color and distance less than $R/2$.
Distinct components are at distance at least $R/4$.
We choose a basis in each $H_C$.
For a singleton, use $k_a$.
For $C=\{a,b\}$ with $\rho(a,b)\leq\tau_0$, use the orthonormal basis $k_a,v_{a,\sigma_a(b)}$ from Lemma~\ref{lem:close-atom-envelope}, assigning center $a$ to both vectors.
If $\rho(a,b)>\tau_0$, use $k_a,k_b$ with their own centers.
In the last case the local Gram matrix has eigenvalues $1\pm\kappa_\beta(a,b)$.
Consequently every local Gram matrix lies between $c_\beta I$ and $2I$, where
\[
c_\beta=1-(1-\tau_0^2)^{\beta/2}>0.
\]
The same bounds include singleton and close-pair blocks.

The list of assigned centers has multiplicity at most two.
Distinct assigned centers in the same two-point component are separated by at least $\operatorname{arctanh}\tau_0$, while centers belonging to distinct components are separated by at least $R/4$.
Thus, after fixing $R\geq R_0$, the distinct assigned centers have a uniform positive separation.
Let $G_0$ be the block-diagonal local Gram matrix and let $E$ be the off-block part of the global Gram matrix.
By Lemmas~\ref{lem:weighted-tail} and \ref{lem:close-atom-envelope},
\[
\sup_i\sum_j|E_{ij}|\frac{w(a_j)}{w(a_i)}
\leq C_\beta(1+R)e^{-(\beta-1)R/4}.
\]
By symmetry, the second weighted Schur quantity is obtained from the same estimate after interchanging the two indices.
The weighted Schur test therefore gives the same bound for $\|E\|_{\ell^2\to\ell^2}$.
Increase $R_0$ so that $\|E\|_{\ell^2\to\ell^2}\leq c_\beta/2$ for every $R\geq R_0$.
The global synthesis map $W$ for the chosen vectors then satisfies
\begin{equation}
\label{eq:global-gram-bounds}
\frac{c_\beta}{2}\|z\|_2^2\leq\|Wz\|_H^2
\leq(2+c_\beta/2)\|z\|_2^2.
\end{equation}
These bounds are first obtained on finite sequences and hence define $W$ on all of $\ell^2$.

Let $S:\mathcal X\to H$ be the summation map $S(x_C)=\sum_Cx_C$.
Comparison with the local Gram matrices and \eqref{eq:global-gram-bounds} gives
\[
\frac{c_\beta}{4}\sum_C\|x_C\|^2
\leq\|Sx\|^2
\leq\frac{2+c_\beta/2}{c_\beta}\sum_C\|x_C\|^2.
\]
The adjoint of $S$ is exactly $\mathcal R$ in \eqref{eq:cluster-restriction}.
Hence $S^*S$ is boundedly invertible, and
\[
\mathcal A=S(S^*S)^{-1},\qquad
\mathcal R\mathcal A=I,\qquad \|\mathcal A\|\leq2c_\beta^{-1/2}.
\]
This proves the required uniformity, including for infinite sequences.
The one-color case is included by taking the second sequence empty.
\end{proof}

\section{Necessary estimates}
\label{sec:necessity}

We now prove the remaining necessary estimates: the strong estimate for $0<p<1$ and the weak estimate for $0<p\leq1$.
Split $\D$ into $\mathcal N=\{z\in\D:\rho_z<\delta\}$ and $\mathcal F=\{z\in\D:\rho_z\geq\delta\}$, where $0<\delta<1$ will be fixed sufficiently small in terms of $\alpha$.
For $\xi\in\{\varphi,\psi\}$, define the restricted cell densities
\[
q_{j,\xi}^{\rm near}
=s_j^2\int_{\mathcal N\cap\xi^{-1}(Q_j)}\rho_z^2\,\dd A_\alpha(z),
\qquad
q_{j,\xi}^{\rm far}
=s_j^2\int_{\mathcal F\cap\xi^{-1}(Q_j)}\rho_z^2\,\dd A_\alpha(z).
\]
Set $q_j^{\rm near}=q_{j,\varphi}^{\rm near}+q_{j,\psi}^{\rm near}$ and $q_j^{\rm far}=q_{j,\varphi}^{\rm far}+q_{j,\psi}^{\rm far}$; then $q_j=q_j^{\rm near}+q_j^{\rm far}$.

\subsection{Nearby image points}

We approximate differences at nearby image points by normalized derivatives.
For $a\in\D$, write $\mathfrak D_af=\Delta(a)f'(a)$.

\begin{lemma}
\label{lem:derivative-interpolation}
For sufficiently large $R=R(\alpha)$, every $R$-separated sequence $\Lambda$ has a uniformly bounded right inverse for
\[
\mathcal S:H\to\ell^2(\Lambda),\qquad
(\mathcal Sf)_a=s_a^{-1}\mathfrak D_af.
\]
Thus there is $\mathcal A_1$ such that $\mathcal S\mathcal A_1=I$ and $\|\mathcal A_1\|\leq A_1(\alpha)$, independently of the sequence.
\end{lemma}

\begin{proof}
Put $L_a=\Delta(a)K_a^{[1]}$ and $e_a^{(1)}=L_a/\|L_a\|$.
By \eqref{eq:derivative-kernel-norm},
\begin{equation}
\label{eq:normalized-derivative-size}
\|L_a\|=c_as_a,\qquad
c_a=\sqrt{\beta(1+\beta|a|^2)},\qquad
\sqrt\beta\leq c_a\leq\sqrt{\beta(1+\beta)}.
\end{equation}
Direct differentiation gives
\[
\langle K_a^{[1]},K_b^{[1]}\rangle
=\frac{\beta(1+\beta\overline a b)}{(1-\overline a b)^{\beta+2}}.
\]
Since $|1+\beta\overline a b|^2\leq(1+\beta|a|^2)(1+\beta|b|^2)$,
\[
|\langle e_a^{(1)},e_b^{(1)}\rangle|
\leq(1-\rho(a,b)^2)^{(\beta+2)/2}
\leq C_\beta e^{-(\beta+2)d(a,b)}.
\]
Disjoint fixed-radius hyperbolic disks show that the number of centers at distance in $[R+j,R+j+1)$ from $a$ is at most $Ce^{2(R+j)}$.
It follows that the sum of the absolute off-diagonal entries in any Gram row is at most $C_\beta e^{-\beta R}$.
Choose $R$ so that this is at most $1/2$.
If $E$ denotes synthesis by the vectors $e_a^{(1)}$ and $G=E^*E$, then $\tfrac12I\leq G\leq\tfrac32I$.
Since the inner product is linear in its first variable,
\[
\langle f,e_a^{(1)}\rangle
=\frac{\mathfrak D_af}{c_as_a}.
\]
Thus, with $C_0=\operatorname{diag}(c_a)$,
\[
\mathcal S=C_0E^*,\qquad
\mathcal A_1=EG^{-1}C_0^{-1}.
\]
The second identity satisfies $\mathcal S\mathcal A_1=I$, and \eqref{eq:normalized-derivative-size} bounds its norm by a constant depending only on $\alpha$.
\end{proof}

For the local expansion below, fix $\delta_0=1/16$ and $r_*=1$, and put
\[
\|f\|_{j,*}^2=\int_{B(a_j,r_*)}|f|^2\,\dd A_\alpha.
\]

\begin{lemma}
\label{lem:near-expansion}
Let $0<\eta\leq\delta_0$, $a\in Q_j$, and $0<\rho(a,b)<\delta\leq\delta_0$.
There are $|\omega(a,b)|=1$ and a linear functional $\mathfrak R_{a,b}$ such that
\begin{align}
f(b)-f(a)&=\omega(a,b)\rho(a,b)\mathfrak D_af
+\rho(a,b)^2\mathfrak R_{a,b}f,
\label{eq:exact-near-expansion}\\
|\mathfrak R_{a,b}f|&\leq Cs_j\|f\|_{j,*},
\label{eq:local-remainder}\\
|\mathfrak D_af-\mathfrak D_{a_j}f|&\leq C\eta s_j\|f\|_{j,*}.
\label{eq:local-freezing}
\end{align}
The constants are independent of $\eta,\delta,a,b,j$ within these fixed ranges.
\end{lemma}

\begin{proof}
Let $w=\sigma_a(b)$, $\omega=-w/|w|$, and $F(\zeta)=s_a^{-1}f(\sigma_a(\zeta))$.
For $0\leq t\leq1$,
\[
d(a_j,\sigma_a(tw))\leq\eta+\operatorname{arctanh}\delta_0<1/4.
\]
The Euclidean segment from $a_j$ to $a$ also lies in $B(a_j,1/4)$, since this hyperbolic disk is a Euclidean disk.
Hyperbolic $1/4$-balls about these points are contained in $B(a_j,r_*)$.
Thus the local derivative estimates \eqref{eq:local-derivatives} along both paths are controlled by $\|f\|_{j,*}$.

Taylor's integral formula gives
\[
f(b)-f(a)=-w\mathfrak D_af
+s_aw^2\int_0^1(1-t)F''(tw)\,\dd t.
\]
This is \eqref{eq:exact-near-expansion} with $\mathfrak R_{a,b}f=s_a\omega^2\int_0^1(1-t)F''(tw)\,\dd t$.
On the fixed disk $|\zeta|\leq\delta_0$,
\[
\sigma_a'(\zeta)=-\frac{\Delta(a)}{(1-\overline a\zeta)^2},
\qquad
\sigma_a''(\zeta)=-\frac{2\overline a\Delta(a)}
{(1-\overline a\zeta)^3},
\]
so both derivatives are $O(\Delta(a))$ there.
Applying \eqref{eq:local-derivatives} with orders one and two, together with \eqref{eq:local-comparability}, gives $|F''(tw)|\leq C\|f\|_{j,*}$ and proves \eqref{eq:local-remainder}.

For \eqref{eq:local-freezing}, parameterize the segment by $u(t)=a_j+t(a-a_j)$.
Then
\[
\frac{\dd}{\dd t}\bigl(\Delta(u(t))f'(u(t))\bigr)
=\frac{\dd}{\dd t}\Delta(u(t))\,f'(u(t))
+\Delta(u(t))f''(u(t))(a-a_j).
\]
The local derivative estimates give
\[
|f'(u(t))|+\Delta(a_j)|f''(u(t))|
\leq Cs_j\Delta(a_j)^{-1}\|f\|_{j,*}.
\]
Moreover, $|\frac{\dd}{\dd t}\Delta(u(t))|\leq2|a-a_j|$ and $|a-a_j|\leq C\eta\Delta(a_j)$.
Integration over $t\in[0,1]$ proves \eqref{eq:local-freezing}.
\end{proof}

The derivative model gives necessary estimates for the near densities.

\begin{proposition}
\label{prop:near-necessity}
There are $\delta>0$ and $\eta_0>0$, depending only on $\alpha$, with the following property.
For every fixed mesh $0<\eta\leq\eta_0$ and every $0<p<2$, the near densities satisfy
\begin{align*}
\|(q_j^{\rm near})^{1/2}\|_{\ell^p}&\leq C\|D\|_{S_p}
&&(D\in S_p),\\
\|(q_j^{\rm near})^{1/2}\|_{\ell^{p,\infty}}&\leq C\|D\|_{S_{p,\infty}}
&&(D\in S_{p,\infty}).
\end{align*}
The constant may depend on $\alpha,p$ and the fixed mesh $\eta$, but not on $\varphi$ or $\psi$.
\end{proposition}

\begin{proof}
Choose $R$, depending only on $\alpha$, large enough for Lemma~\ref{lem:derivative-interpolation} and so that the $r_*$-balls about points of an $R$-separated sequence are disjoint.
Partition the lattice into finitely many $R$-separated color classes.
Work initially with $0<\delta,\eta\leq\delta_0$.
For $\varphi$, put
\[
\Omega_{j,\varphi}=\mathcal N\cap\varphi^{-1}(Q_j),\qquad
\tau_j^2=\int_{\Omega_{j,\varphi}}\rho_z^2\,\dd A_\alpha,
\qquad \lambda_j=s_j\tau_j.
\]
In particular, $\lambda_j^2=q_{j,\varphi}^{\rm near}$.
For $\tau_j>0$, define
\begin{equation}
\label{eq:near-output-compression}
(Q_\varphi g)_j=-\tau_j^{-1}
\int_{\Omega_{j,\varphi}}\overline{\omega(\varphi(z),\psi(z))}
\rho_z g(z)\,\dd A_\alpha(z).
\end{equation}
Rows with $\tau_j=0$ are set equal to zero.
On $\{\rho_z=0\}$, set the phase to one and the remainder to zero.
Cauchy--Schwarz and disjointness give $\|Q_\varphi\|\leq1$.

Apply \eqref{eq:exact-near-expansion} with $a=\varphi(z)$ and $b=\psi(z)$.
Since $\DP f=f(a)-f(b)$, the phase in \eqref{eq:near-output-compression} cancels the principal term.
For $\tau_j>0$, the resulting error is
\begin{align*}
E_jf={}&\tau_j^{-1}\int_{\Omega_{j,\varphi}}\rho_z^2
\bigl(\mathfrak D_{\varphi(z)}f-\mathfrak D_{a_j}f\bigr)
\,\dd A_\alpha(z)\\
&+\tau_j^{-1}\int_{\Omega_{j,\varphi}}
\overline{\omega(\varphi(z),\psi(z))}\rho_z^3
\mathfrak R_{\varphi(z),\psi(z)}f\,\dd A_\alpha(z).
\end{align*}
Set $E_j=0$ if $\tau_j=0$.
Using \eqref{eq:local-remainder}--\eqref{eq:local-freezing} and $\int_{\Omega_{j,\varphi}}\rho_z^3\,\dd A_\alpha\leq\delta\tau_j^2$, we obtain
\begin{equation}
\label{eq:near-row-decomposition}
(Q_\varphi\DP f)_j=\tau_j\mathfrak D_{a_j}f+E_jf,
\qquad
|E_jf|\leq C(\eta+\delta)\lambda_j\|f\|_{j,*}.
\end{equation}

On a color class $\Lambda_\nu$, divide $E_j$ by $\lambda_j$ when positive and use zero otherwise.
This defines an operator $\mathcal E_{\nu,\varphi}:H\to\ell^2(\Lambda_\nu)$ with
\begin{equation}
\label{eq:near-error-norm}
\|\mathcal E_{\nu,\varphi}\|\leq C(\eta+\delta),
\end{equation}
because the $r_*$-balls are disjoint on that color.
Let $T_{\nu,\varphi}$ be the corresponding rows of $Q_\varphi\DP$.
Its singular values are bounded by those of $D$ whenever $D$ is compact.

Choose $\delta,\eta_0\leq\delta_0$ so that $CA_1(\eta_0+\delta)\leq1/2$ in \eqref{eq:near-error-norm}.
For a finite subset $F$ of the color class, let $I_F$ be coordinate inclusion, $P_F=I_F^*$, $A_F=\mathcal A_1I_F$, and $\Lambda_F=\operatorname{diag}(\lambda_j:j\in F)$.
Equation \eqref{eq:near-row-decomposition} and $\mathcal S\mathcal A_1=I$ give
\[
P_FT_{\nu,\varphi}A_F=\Lambda_F(I+B_F),\qquad
B_F=P_F\mathcal E_{\nu,\varphi}\mathcal A_1I_F,
\qquad \|B_F\|\leq1/2.
\]
Inverting $I+B_F$ yields
\[
s_n(\Lambda_F)\leq2A_1s_n(D),\qquad n\geq1.
\]
Summation gives the strong estimate on $F$, and the definition of the weak quasi-norm gives the weak estimate.
Increasing finite sets proves both estimates on the color class.
Repeat the argument with $\psi$ in place of $\varphi$, using $-\DP$.
Finally combine the finitely many colors and the two density sequences.
For strong estimates use $(a+b)^{p/2}\leq a^{p/2}+b^{p/2}$, and for weak estimates use the finite level-set union bound.
The number of colors depends only on $R$ and the fixed mesh, so all constants remain independent of the symbols.
\end{proof}
	
\subsection{Finite edge blocks}

We retain the cutoff $\delta$ and turn to the complementary set $\mathcal F$.
We construct finite matrices that control the far densities and can be compared with restrictions of $\DP$.

Choose $R=R(\alpha)$ large enough for Proposition~\ref{prop:cluster-interpolation} and Lemma~\ref{lem:derivative-interpolation}, and take $R\geq4$.
For a mesh $0<\eta\leq\tfrac14\operatorname{arctanh}\delta$, partition the lattice into finitely many $R$-separated color classes $\Lambda_\nu$.
Define
\[
\Omega_{jk}=\mathcal F\cap\varphi^{-1}(Q_j)\cap\psi^{-1}(Q_k),
\qquad
\gamma_{jk}=\int_{\Omega_{jk}}\dd A_\alpha(z).
\]
Retain the ordered edges $e=(j,k)$ for which $\gamma_{jk}>0$, and write $\gamma_e=\gamma_{jk}$ and $\Omega_e=\Omega_{jk}$.
Group these edges by the unordered pair of their endpoint colors, retaining the orientation of each edge.
Let $\mathcal E_{\nu\kappa}$ denote the block with color pair $\{\nu,\kappa\}$, where $\nu\leq\kappa$.

Fix one block and use the clusters $C$ from Proposition~\ref{prop:cluster-interpolation} for $\Gamma=\Lambda_\nu\cup\Lambda_\kappa$; if $\nu=\kappa$, use that single sequence.
For a finite edge set $F\subset\mathcal E_{\nu\kappa}$, let $\mathcal X_F$ be the direct sum of all cluster subspaces $H_C$ touched by $F$.
Each such subspace is included in full, including the second kernel point of a two-point cluster even if that point is not an endpoint of $F$.
Let $C(j)$ be the cluster containing $a_j$ and define $M_F:\mathcal X_F\to\ell^2(F)$ by
\begin{equation}
\label{eq:far-principal-matrix}
(M_Fx)_{(j,k)}=\sqrt{\gamma_{jk}}
\bigl(\langle x_{C(j)},K_{a_j}\rangle-
\langle x_{C(k)},K_{a_k}\rangle\bigr).
\end{equation}
Put
\[
h_{j,F}=s_j^2\sum_{\substack{e\in F\\e\ni j}}\gamma_e,
\]
with zero assigned to unused vertices.

\begin{lemma}
\label{lem:far-principal-lower}
Uniformly in $F$, one has
\begin{align}
\sum_jh_{j,F}^{p/2}&\leq C_{\alpha,p,\delta}\|M_F\|_{S_p}^p,
&&0<p<2,\label{eq:far-principal-strong}\\
\|(h_{j,F}^{1/2})\|_{\ell^{p,\infty}}
&\leq C_{\alpha,p,\delta}\|M_F\|_{S_{p,\infty}},
&&p>0.\label{eq:far-principal-weak}
\end{align}
\end{lemma}

\begin{proof}
For every retained edge $e=(j,k)$, the triangle inequality and the mesh restriction give
\begin{equation}
\label{eq:far-center-separation}
d(a_j,a_k)\geq\operatorname{arctanh}\delta-2\eta
\geq\tfrac12\operatorname{arctanh}\delta,\qquad
\rho(a_j,a_k)\geq\delta_1:=\tanh(\tfrac12\operatorname{arctanh}\delta)>0.
\end{equation}
In particular, there are no loops.
Choose an orthonormal basis of each $H_C$.
Each row of $M_F$ meets at most two cluster subspaces, each of dimension at most two, so it has at most four nonzero entries.
Let $d_{C,l}$ be the squared column norms in these bases and set $e_C=\sum_l d_{C,l}$.

For an edge joining two different clusters, its contribution to the cluster of $a_j$ is exactly $\gamma_e\|K_{a_j}\|^2$.
For an edge inside $C=\{a_j,a_k\}$, its entire contribution to $e_C$ is
\[
\gamma_e\|K_{a_j}-K_{a_k}\|^2
\asymp_\alpha\gamma_e\rho(a_j,a_k)^2(s_j^2+s_k^2).
\]
By \eqref{eq:far-center-separation} and Lemma~\ref{lem:kernel-difference}, in both cases
\begin{equation}
\label{eq:cluster-energy-comparison}
c_{\alpha,\delta}\sum_{j\in C}h_{j,F}
\leq e_C\leq C_\alpha\sum_{j\in C}h_{j,F}.
\end{equation}
For $\theta=p/2<1$, each cluster has at most two vertices and two coordinates, so
\[
\sum_{C,l}d_{C,l}^{\theta}
\geq\sum_Ce_C^\theta
\geq c_{\alpha,p,\delta}\sum_jh_{j,F}^\theta.
\]
Lemma~\ref{lem:sparse-strong} with row bound four proves the strong estimate.

For the weak estimate, if $h_{j,F}^{1/2}>t$, then \eqref{eq:cluster-energy-comparison} forces at least one coordinate of its cluster to have $d_{C,l}^{1/2}>ct$.
Each cluster accounts for at most two vertices.
Hence
\[
\|(h_{j,F}^{1/2})\|_{\ell^{p,\infty}}
\leq C\|(d_{C,l}^{1/2})\|_{\ell^{p,\infty}}.
\]
Apply Lemma~\ref{lem:sparse-weak}.
\end{proof}
	
\subsection{Taylor estimates and absorption}

We compare $M_F$ with a restriction of $\DP$ by expanding around the lattice centers.
The Taylor remainder can be absorbed uniformly over finite edge sets when the mesh is sufficiently small.

\begin{proposition}
\label{prop:far-necessity}
Fix $0<\delta<1$ and $0<p<2$.
There exists $\eta_0=\eta_0(\alpha,p,\delta)>0$ such that, for every fixed mesh $0<\eta\leq\eta_0$,
\[
\|(q_j^{\rm far})^{1/2}\|_{\ell^{p,\infty}}
\leq C\|D\|_{S_{p,\infty}},\qquad D\in S_{p,\infty}.
\]
If in addition $0<p<1$, then
\[
\|(q_j^{\rm far})^{1/2}\|_{\ell^p}\leq C\|D\|_{S_p},
\qquad D\in S_p.
\]
Here $C$ may depend on $\alpha,p,\delta$ and the fixed mesh $\eta$, but not on $\varphi$ or $\psi$.
\end{proposition}

\begin{proof}
Fix $r_0=1/4$ and initially take $\eta<r_0$ with $\eta\leq\tfrac14\operatorname{arctanh}\delta$.
Work in one color-pair block $\mathcal E_{\nu\kappa}$ and let $F$ be a finite subset of its edges.
Compose the right inverse from Proposition~\ref{prop:cluster-interpolation} with the isometric zero extension $\mathcal X_F\hookrightarrow\mathcal X$ to obtain $A_F:\mathcal X_F\to H$.
Then $\|A_F\|\leq A_0$ and $P_CA_Fx=x_C$ for every cluster occurring in $\mathcal X_F$.
Define
\[
(Q_Fg)_e=\gamma_e^{-1/2}\int_{\Omega_e}g\,\dd A_\alpha,
\qquad T_F=Q_F\DP.
\]
The sets $\Omega_e$ are disjoint, so Cauchy--Schwarz gives $\|Q_F\|\leq1$ and hence $s_n(T_F)\leq s_n(D)$ whenever $D$ is compact.

Expand
\[
s_j^{-1}f(\sigma_{a_j}(\zeta))
=\sum_{l\geq0}c_{j,l}(f)\zeta^l,\qquad |\zeta|<1.
\]
For all endpoints of $F$, let $C_lf=(c_{j,l}(f))_j$.
Parseval on $|\zeta|<r_0$, followed by the automorphism change of variables, gives
\[
|c_{j,l}(f)|^2\leq C_{\alpha,r_0}(l+1)r_0^{-2l}
\int_{B(a_j,\operatorname{arctanh}r_0)}|f|^2\,\dd A_\alpha.
\]
The disks are disjoint within either color and have multiplicity at most two in their union.
Thus
\begin{equation}
\label{eq:coefficient-analysis-bound}
\|C_l\|\leq C_{\alpha,r_0}\sqrt{l+1}\,r_0^{-l},
\end{equation}
uniformly in $F$ and the mesh.

For $e=(j,k)$, put
\[
\theta_{e,j,l}=\gamma_e^{-1}\int_{\Omega_e}
\sigma_{a_j}(\varphi(z))^l\,\dd A_\alpha(z),\qquad
\theta_{e,k,l}=\gamma_e^{-1}\int_{\Omega_e}
\sigma_{a_k}(\psi(z))^l\,\dd A_\alpha(z).
\]
Both moduli are at most $\eta^l$, since $\tanh\eta\leq\eta$.
Let $B_{F,l}$ be the matrix on all endpoint coordinates of $F$ with row
\[
(B_{F,l}y)_e=\sqrt{\gamma_e}
(s_j\theta_{e,j,l}y_j-s_k\theta_{e,k,l}y_k).
\]
It has at most two nonzero entries in each row and column energies at most $\eta^{2l}h_{j,F}$.
The coefficient expansion gives
\begin{equation}
\label{eq:exact-far-transfer}
T_FA_F=M_F+R_F,\qquad
R_F=\sum_{l\geq1}B_{F,l}C_lA_F.
\end{equation}
Indeed, for every endpoint $a_j$,
\[
(A_Fx)(a_j)=\langle P_{C(j)}A_Fx,K_{a_j}\rangle
=\langle x_{C(j)},K_{a_j}\rangle,
\]
so the term with $l=0$ is exactly \eqref{eq:far-principal-matrix}.
For the remainder, \eqref{eq:sparse-pointwise-upper} gives $\|B_{F,l}\|\leq\sqrt2\eta^l\sup_jh_{j,F}^{1/2}$.
Together with \eqref{eq:coefficient-analysis-bound} and $\eta<r_0$, this proves operator-norm convergence in \eqref{eq:exact-far-transfer}.

For the strong estimate, suppose that $0<p<1$ and $D\in S_p$, and write $E_F=\sum_jh_{j,F}^{p/2}$.
Lemma~\ref{lem:sparse-strong} and the ideal property give
\[
\|B_{F,l}C_lA_F\|_{S_p}^p
\leq C_{\alpha,r_0}^pA_0^p(l+1)^{p/2}
(\eta/r_0)^{lp}E_F.
\]
Using the $p$-subadditivity of the Schatten quasi-norm and the convergence just proved, we obtain
\begin{equation}
\label{eq:far-remainder-strong}
\|R_F\|_{S_p}^p
\leq C_{\alpha,r_0}^pA_0^p\Phi_p(\eta/r_0)E_F,\qquad
\Phi_p(t)=\sum_{l\geq1}(l+1)^{p/2}t^{lp}.
\end{equation}
Equations~\eqref{eq:far-principal-strong}, \eqref{eq:exact-far-transfer}, and \eqref{eq:far-remainder-strong} imply
\[
E_F\leq C_{\alpha,p,\delta}
\left(A_0^p\|D\|_{S_p}^p+
C_{\alpha,r_0}^pA_0^p\Phi_p(\eta/r_0)E_F\right).
\]
Since $\Phi_p(t)\to0$ as $t\downarrow0$, a sufficiently small mesh makes the coefficient of $E_F$ on the right at most $1/2$.
It follows that $E_F\leq C\|D\|_{S_p}^p$, uniformly in $F$ and the symbols.

For the weak estimate, suppose that $D\in S_{p,\infty}$ and put $W_F=\|(h_{j,F}^{1/2})\|_{\ell^{p,\infty}}$, which is finite.
Equation \eqref{eq:sparse-pointwise-upper} and the ideal property give
\[
\|B_{F,l}C_lA_F\|_{S_{p,\infty}}
\leq C_{\alpha,r_0}A_0\sqrt{l+1}(\eta/r_0)^lW_F.
\]
Put $K_p=2^{1/p}$.
The singular-value inequality $s_{2n-1}(X+Y)\leq s_n(X)+s_n(Y)$ gives
\begin{equation}
\label{eq:weak-quasitriangle}
\|X+Y\|_{S_{p,\infty}}
\leq K_p\bigl(\|X\|_{S_{p,\infty}}+\|Y\|_{S_{p,\infty}}\bigr);
\end{equation}
see \cite{Simon2005}.
Finite backward iteration yields
\[
\left\|\sum_{l=1}^N B_{F,l}C_lA_F\right\|_{S_{p,\infty}}
\leq C_{\alpha,r_0}A_0W_F
\sum_{l=1}^N K_p^l\sqrt{l+1}(\eta/r_0)^l.
\]
If $K_p\eta<r_0$, operator-norm convergence and continuity of each singular value allow $N\to\infty$.
Consequently,
\begin{equation}
\label{eq:far-remainder-weak}
\|R_F\|_{S_{p,\infty}}
\leq C_{\alpha,r_0}A_0\Psi_p(\eta/r_0)W_F,\qquad
\Psi_p(t)=\sum_{l\geq1}K_p^l\sqrt{l+1}\,t^l.
\end{equation}
By \eqref{eq:far-principal-weak}, \eqref{eq:weak-quasitriangle}, and \eqref{eq:far-remainder-weak},
\[
W_F\leq C_{\alpha,p,\delta}K_p
\left(A_0\|D\|_{S_{p,\infty}}
+C_{\alpha,r_0}A_0\Psi_p(\eta/r_0)W_F\right).
\]
Since $\Psi_p(t)\to0$ as $t\downarrow0$, decreasing the mesh absorbs the last term and gives $W_F\leq C\|D\|_{S_{p,\infty}}$ uniformly in $F$.
Choose $\eta_0$ so that $\eta_0\leq\tfrac14\operatorname{arctanh}\delta$, $K_p\eta_0<r_0$, and the coefficient of $W_F$ on the right is at most $1/2$ for every $0<\eta\leq\eta_0$.
When $p<1$, also require that the coefficient of $E_F$ in the strong estimate is at most $1/2$.
All these choices are independent of $F$ and the symbols.

As $F$ increases to exhaust the fixed block,
\[
h_{j,F}\uparrow h_j^{\nu\kappa}
:=s_j^2\sum_{\substack{e\in\mathcal E_{\nu\kappa}\\e\ni j}}\gamma_e.
\]
Monotone convergence proves the strong estimate on the block.
For the weak estimate, let $J$ be any finite subset of $\{j:(h_j^{\nu\kappa})^{1/2}>t\}$.
There is a sufficiently large finite edge set $F$ for which $h_{j,F}^{1/2}>t$ for all $j\in J$.
Hence $t|J|^{1/p}\leq W_F\leq C\|D\|_{S_{p,\infty}}$.
Taking the supremum over finite $J$ and $t>0$ gives the weak estimate for the block.

The finitely many block degrees add to
\[
h_j=\sum_{\nu\leq\kappa}h_j^{\nu\kappa}
=s_j^2\sum_{e\ni j}\gamma_e.
\]
Combine the block estimates using subadditivity of $t^{p/2}$ in the strong case and a finite level-set union bound in the weak case.
The number of blocks enters only these final comparisons.
Since $\delta^2\leq\rho_z^2\leq1$ on $\mathcal F$,
\[
\delta^2h_j\leq q_j^{\rm far}\leq h_j.
\]
This proves both assertions.
\end{proof}
	
	\Needspace{9\baselineskip}
\section{Toeplitz estimates and sufficiency}
\label{sec:toeplitz}

We prove both sufficient directions by dominating $D^*D$ with a positive Toeplitz operator.
We first record the spectral estimates and then construct the majorant used in both proofs.

For a positive compact operator $T$, let $\lambda_n(T)$ denote its eigenvalues in nonincreasing order, counted with multiplicity.
For a nonnegative locally integrable function $g$, let $\Toep_g$ be the positive Toeplitz operator, when bounded, with form
\[
\langle\Toep_gf,f\rangle=\int_\D g(u)|f(u)|^2\,\dd A_\alpha(u).
\]
Choose a fixed sufficiently fine lattice of Euclidean disks $B_i=\{z\in\D:|z-z_i|<c\Delta(z_i)\}$, with a fixed small $c>0$, of the type used in \cite[Definition~2.5]{ElFallahElIbbaoui2022}, and set
\[
A(B_i)=\int_{B_i}\dd A,\qquad
b_i(g)=\frac1{A(B_i)}\int_{B_i}g\,\dd A.
\]

\begin{lemma}
\label{lem:toeplitz-sufficiency}
For every $s>0$, if $(b_i(g))\in\ell^s$, then
\begin{equation}
\|\Toep_g\|_{S_s}\leq C_{\alpha,s}\|(b_i(g))\|_{\ell^s}.
\label{eq:toeplitz-strong}
\end{equation}
Moreover, if $(b_i(g))\in\ell^{s,\infty}$, then
\begin{equation}
\|\Toep_g\|_{S_{s,\infty}}\leq
C_{\alpha,s}\|(b_i(g))\|_{\ell^{s,\infty}}.
\label{eq:toeplitz-weak}
\end{equation}
The constants may also depend on the fixed lattice used to define the averages $b_i(g)$, but not on $g$.
\end{lemma}

\begin{proof}
Either sequence assumption implies that $\sup_i b_i(g)<\infty$ and that $b_i(g)\to0$ outside compact subsets of $\D$.
Local submean estimates on the lattice disks and bounded overlap give
\[
\int_\D g|f|^2\,\dd A_\alpha
\leq C\bigl(\sup_i b_i(g)\bigr)\|f\|_H^2.
\]
Applying the same estimate to the union of disks outside a fixed compact set shows that the boundary-tail form norm tends to zero.
On a compact subset, bounded sets in $A_\alpha^2$ are relatively compact for local uniform convergence, and $g$ is integrable there.
Hence the corresponding form is compact, and therefore $\Toep_g$ is compact.

The strong Schatten criterion is classical \cite{Luecking1987}.
For quantitative estimates, we use \cite[Theorem~4.5]{ElFallahElIbbaoui2022}.
In their convention, the squared Bergman weight is $(\alpha+1)\Delta(z)^\alpha/\pi$ and the Toeplitz measure is $\dd\nu(z)=g(z)\,\dd x\,\dd y$, so their Toeplitz form agrees with the form defining $\Toep_g$.
Replacing Euclidean area by $\dd A=\pi^{-1}\dd x\,\dd y$ multiplies the numerator and denominator of each local average by the same constant.
Thus $b_i(g)$ are the geometric averages in their estimates.

For every $\theta>0$, apply that theorem with $h(t)=t^\theta$.
This function is increasing, satisfies $h(0)=0$, and $h(t^{1/\theta})=t$ is convex, as required.
The upper estimate gives
\[
\sum_{i=1}^n\lambda_i(\Toep_g)^\theta
\leq C_{\alpha,\theta}\sum_{i=1}^n b_i^*(g)^\theta,
\qquad n\geq1,
\]
where the constant also depends on the fixed lattice but not on $g$.
Taking $\theta=s$ and then $n\to\infty$ proves \eqref{eq:toeplitz-strong}.

For the weak estimate, let $M=\|(b_i(g))\|_{\ell^{s,\infty}}$ and choose $0<\theta<\min\{s,1\}$.
Then $b_i^*(g)\leq Mi^{-1/s}$ and
\[
n\lambda_n(\Toep_g)^\theta
\leq\sum_{i=1}^n\lambda_i(\Toep_g)^\theta
\leq C M^\theta\sum_{i=1}^n i^{-\theta/s}
\leq C_{\alpha,s,\theta}M^\theta n^{1-\theta/s}.
\]
Thus $\lambda_n(\Toep_g)\leq CMn^{-1/s}$, proving \eqref{eq:toeplitz-weak}.
\end{proof}

To apply these estimates, we construct a Toeplitz majorant whose local averages are controlled by the cell densities.

\begin{lemma}
\label{lem:toeplitz-majorant}
There is a fixed radius $r_2>0$ such that
\[
g(u)=\int_\D\mathbf1_{B(a,r_2)}(u)\Delta(a)^{-\beta}\,\dd\mu(a)
\]
is nonnegative and locally integrable, and
\begin{equation}
\label{eq:positive-domination}
D^*D\leq C\Toep_g
\end{equation}
whenever the defining form of $\Toep_g$ is bounded.
For every $p>0$ and every fixed cell lattice,
\begin{align}
\|(b_i(g))\|_{\ell^{p/2}}&\leq C\|(q_j)\|_{\ell^{p/2}},
\label{eq:majorant-strong-sequence}\\
\|(b_i(g))\|_{\ell^{p/2,\infty}}&\leq
C\|(q_j)\|_{\ell^{p/2,\infty}}.
\label{eq:majorant-weak-sequence}
\end{align}
The constants depend only on $\alpha$, $p$, and the fixed lattices, and are independent of $\varphi$ and $\psi$.
\end{lemma}

\begin{proof}
There are fixed radii $0<r_1<r_2$ for which
\begin{equation}
\label{eq:local-difference-majorization}
|f(a)-f(b)|^2\leq C\rho(a,b)^2
\left(s_a^2\int_{B(a,r_2)}|f|^2\,\dd A_\alpha
+s_b^2\int_{B(b,r_2)}|f|^2\,\dd A_\alpha\right).
\end{equation}
Indeed, when $d(a,b)\geq r_1$, the factor $\rho(a,b)$ is bounded below, and \eqref{eq:local-derivatives} with $m=0$ gives the estimate.
When $d(a,b)<r_1$, integrate $f'$ along the segment from $a$ to $b$ and apply \eqref{eq:local-derivatives} with $m=1$ on a larger fixed disk.
The two depths are comparable by \eqref{eq:local-comparability}, and $|a-b|\lesssim\rho(a,b)\Delta(a)$.
This proves \eqref{eq:local-difference-majorization}; compare \cite[Section~4]{DuTongYang2026}.

For $u$ in a compact subset of $\D$, the points $a$ contributing to $g(u)$ stay in a larger compact subset.
Since $\mu$ is finite, $g$ is locally bounded and hence locally integrable.
Tonelli's theorem and \eqref{eq:local-difference-majorization} yield \eqref{eq:positive-domination} whenever the Toeplitz form is bounded.
This majorization is also used for the strong sufficient condition in \cite{ChoeChoiKooPark2026}; see also \cite{DuTongYang2026,HuYang2026}.

By \eqref{eq:local-comparability},
\begin{equation}
\label{eq:majorant-average}
g(u)\leq C\Delta(u)^{-\beta}\mu(B(u,r_2)).
\end{equation}
For each disk $B_i$, fixed-radius geometry gives a set $\mathcal J_i$ of nearby cell indices, with uniformly bounded cardinality and reverse incidence multiplicity.
The corresponding centers have depths comparable to $\Delta(z_i)$.
Equation \eqref{eq:majorant-average} therefore implies
\[
b_i(g)\leq C\sum_{j\in\mathcal J_i}q_j.
\]
The finite sum and level-set arguments in Lemma~\ref{lem:discretization} prove \eqref{eq:majorant-strong-sequence} and \eqref{eq:majorant-weak-sequence}.
\end{proof}
	
We can now complete the strong Schatten characterization.

\begin{proof}[Proof of Theorem~\ref{thm:strong}]
For necessity when $1\leq p<2$, use Proposition~\ref{prop:banach-necessity}.
For $0<p<1$, choose the cutoff and the mesh as in Propositions~\ref{prop:near-necessity} and \ref{prop:far-necessity}.
The near and far measures add to $\mu$, so
\[
\sum_jq_j^{p/2}\leq
\sum_j(q_j^{\rm near})^{p/2}+
\sum_j(q_j^{\rm far})^{p/2}
\leq C\|D\|_{S_p}^p.
\]
Lemma~\ref{lem:discretization} transfers the estimate to the lattice in the statement.

Conversely, suppose that $\sum_jm_j^{p/2}<\infty$ and let $g$ be given by Lemma~\ref{lem:toeplitz-majorant}.
Then \eqref{eq:strong-discretization}, \eqref{eq:majorant-strong-sequence}, and Lemma~\ref{lem:toeplitz-sufficiency} imply $\Toep_g\in S_{p/2}$ and
\[
\|\Toep_g\|_{S_{p/2}}\leq
C\|(q_j)\|_{\ell^{p/2}}.
\]
Since $\Toep_g$ is compact, the positive domination \eqref{eq:positive-domination} implies that $D$ is compact.
The min--max principle then gives
\[
s_n(D)^2=\lambda_n(D^*D)\leq C\lambda_n(\Toep_g),\qquad n\geq1.
\]
Consequently,
\[
\|D\|_{S_p}^2=\|D^*D\|_{S_{p/2}}
\leq C\|\Toep_g\|_{S_{p/2}}.
\]
Together with $\|(q_j)\|_{\ell^{p/2}}=\|(q_j^{1/2})\|_{\ell^p}^2$, this proves sufficiency and the other half of \eqref{eq:main-strong-norm}.

The sufficiency argument applies to every $p>0$, and Proposition~\ref{prop:banach-necessity} gives necessity for $p\geq2$ as well.
Hence the same equivalence and two-sided norm estimate hold for every $p>0$, also recovering the previously known $p\geq2$ characterization in \cite{ChoeChoiKooPark2026}.
\end{proof}

\begin{proof}[Proof of Theorem~\ref{thm:weak}]
For necessity when $0<p\leq1$, Propositions~\ref{prop:near-necessity} and \ref{prop:far-necessity} give the required weak bounds for the two restricted density sequences.
If $q_j^{1/2}>t$, then either $(q_j^{\rm near})^{1/2}>t/\sqrt2$ or $(q_j^{\rm far})^{1/2}>t/\sqrt2$.
A level-set union bound followed by \eqref{eq:weak-discretization} proves necessity in this range.

For $p>1$, necessity follows from \eqref{eq:weak-necessity-large-p} and \eqref{eq:weak-discretization}.

For sufficiency, let $p>0$ and let $g$ be given by Lemma~\ref{lem:toeplitz-majorant}.
If $(m_j^{1/2})\in\ell^{p,\infty}$, then $(q_j)\in\ell^{p/2,\infty}$.
Equations \eqref{eq:majorant-weak-sequence} and \eqref{eq:toeplitz-weak} give $\Toep_g\in S_{p/2,\infty}$.
Since $\Toep_g$ is compact, \eqref{eq:positive-domination} implies that $D$ is compact.
The min--max principle then gives
\[
s_n(D)^2\leq C\lambda_n(\Toep_g),
\]
so $D\in S_{p,\infty}$ and
\[
\|D\|_{S_{p,\infty}}^2
=\|D^*D\|_{S_{p/2,\infty}}
\leq C\|\Toep_g\|_{S_{p/2,\infty}}
\leq C\|(q_j^{1/2})\|_{\ell^{p,\infty}}^2.
\]
This proves \eqref{eq:main-weak-norm}.
\end{proof}

The argument includes degenerate symbols.
If $\varphi=\psi$, all quantities vanish.
If both symbols are constant, the difference has rank at most one and the joint measure is supported on at most two points.
Since $\mu$ includes both pullback measures, the criterion also covers cases in which only one symbol approaches the boundary.

\declarationheading{Author contributions}
Both authors contributed equally to the research and writing of this manuscript.
Both authors read and approved the final manuscript.

\declarationheading{Conflict of interest}
The authors declare that they have no conflicts of interest.

\declarationheading{Data availability statement}
No datasets were generated or analyzed during the current study.

\declarationheading{Declaration of Generative AI Use}
The authors used artificial-intelligence tools for language editing,
\LaTeX\ formatting, and assistance with local mathematical
reasoning.
All mathematical content was reviewed and verified by the authors,
who take full responsibility for the final manuscript.

\end{document}